\documentclass[a4paper,12pt,titlepage,french]{smfart}

\usepackage{vmargin}
\usepackage{amssymb, amsfonts}
\usepackage{amsmath}
\usepackage[frenchb]{babel}
\usepackage[latin1]{inputenc}
\usepackage[ec]{aeguill}
\usepackage[pdftex]{graphics}
\usepackage{xy}
\usepackage{xypic}
\usepackage{graphicx}
\usepackage{color}
\usepackage{pstricks}

\newcommand{\field}[1]{\mathbb{#1}}
\newcommand{\rrr}{\field{R}}
\newcommand{\zz}{\frac{\zzz}{2\zzz}}
\newcommand{\ccc}{\field{C}}
\newcommand{\nnn}{\field{N}}
\newcommand{\zzz}{\field{Z}}

\newcommand{\qqq}{\field{Q}}

\newcommand{\A}{\mathcal{A}}

\renewcommand{\O}{\mathcal{O}}
\renewcommand{\L}{\mathcal{L}}

\newcommand{\XR}{X_{\Delta}(\field{R})}

\newcommand{\XC}{X_{\Delta}(\field{C})}

\newcommand{\Vsig}{\displaystyle \frac{V}{[\sigma]_2}}
\newcommand{\Vtau}{\displaystyle \frac{V}{[\tau]_2}}

\newcommand{\BB}{\textstyle \text{{\large $\bigwedge^*$}}}
\newcommand{\Bq}{\textstyle \text{{\large $\bigwedge^q$}}}

\newcommand{\dimm}{\text{dim}}

\newcommand{\Vect}{\text{Vect}}
\newcommand{\card}{\text{Card}}

\newcommand{\Iq}{\displaystyle \frac{I^q}{I^{q+1}}}
\newcommand{\I}{\displaystyle \frac{I}{I^2}}

\newtheorem{theoreme}{Th{\'e}or{\`e}me}[section]
\newtheorem{proposition}[theoreme] {Proposition}
\newtheorem{lemme}[theoreme]{Lemme}
\newtheorem{corollaire}[theoreme]{Corollaire}
\newtheorem{definition}[theoreme]{D{\'e}finition}

\newenvironment{remarque}{\textbf{Remarque:\\}}{}

\newenvironment{preuve}{\textbf{D\'emonstration:\\}}{}

\begin{document}

\title{Maximalit\'e des vari\'et\'es toriques de dimension 4}
\date {Angers}
\author{Alexandre Sine}

\address{D\'epartement de Math\'ematiques, UMR CNRS 6093,
Universit\'e d'Angers, 2 Boulevard Lavoisier, 49045 Angers Cedex
01, France.} \email{sine@tonton.univ-angers.fr} \urladdr{}

\begin{abstract}
Une vari\'et\'e alg\'ebrique complexe d\'efinie sur les r\'eels
est dite maximale si la somme de ses nombres de Betti pour
l'homologie de Borel Moore \`a coefficients dans $\zz$ est \'egale
\`a la somme de nombres de Betti de sa partie r\'eelle. On
montrera dans cet article
que les vari\'et\'es toriques de dimension 4 sont maximales.\\

{\sc Mots-cl\'es:} vari\'et\'e torique r\'eelle, homologie,
maximalit\'e

\end{abstract}

\begin{altabstract}
A complex algebraic variety defined over the reals is maximal when
the sum of its Betti numbers for Borel Moore homology with $\zz$
coefficients coincides with the sum of the Betti numbers of its
real part. We will show in this paper that toric
varieties of dimension 4 are maximal.\\

{\sc Index words:} real toric variety, homology, maximality
\end{altabstract}

\maketitle

\section{Introduction}

Dans cet article, on va s'int\'eresser \`a la question suivante:\\
\hglue1truecm{\emph{Est ce que les vari\'et\'es toriques sont maximales?}}\\
Ici, une vari\'et\'e alg\'ebrique complexe $X:=X(\ccc)$ d\'efinie
sur $\rrr$ est maximale lorsque, pour l'homologie de Borel Moore
\`a coefficients dans $\zz$, la somme de ses nombres de Betti
co\"incide avec la somme correspondante pour sa partie r\'eelle,
$X(\rrr)$. Ce qui peut \^etre r\'esum\'e par l'\'egalit\'e:

$$\sum_i\dimm ~H^{BM}_i(X(\rrr),\zz) = \sum_i\dimm ~
H^{BM}_i(X(\ccc),\zz)$$

Dans \cite{Mc} (page 11), les 4 auteurs ont prouv\'e que ce
r\'esultat est vrai jusqu'en dimension 3. Ceci f\^ut obtenu en
construisant une suite spectrale $G^*_{pq}(X(\rrr))$ convergeant
vers $H^{BM}_*(X(\rrr),\zz)$, en la comparant \`a une autre suite
spectrale $E^*_{pq}(X(\ccc))$ convergeant, elle, vers
$H^{BM}_*(X(\ccc),\zz)$ et en utilisant l'in\'egalit\'e de
Thom-Smith:

$$\sum_i\dimm~ H^{BM}_i(X(\rrr),\zz) \leq \sum_i\dimm~
H^{BM}_i(X(\ccc),\zz)$$

Le point cl\'e \'etant la d\'eg\'en\'erescence \`a l'ordre 1 de la
suite spectrale $G^*_{pq}(X(\rrr))$ (ie:
$G^{\infty}_{pq}(X(\rrr))=G^1_{pq}(X(\rrr))$).\\

Dans \cite{Ho}, V.Hower a montr\'e que cette suite d\'eg\'en\`ere
\`a l'ordre 1 lorsque $X$ est une vari\'et\'e torique projective
dont le polytope $\Delta$ lui \'etant associ\'e v\'erifie:
\begin{itemize}
    \item $\Delta$ est r\'eflexif.
    \item L'\'eventail normal du polytope dual de $\Delta$ est r\'egulier modulo 2.
\end{itemize}
Dans \cite{Ho2}, elle a donn\'e un exemple de vari\'et\'e torique
de dimension 6 n'\'etant pas maximal. Cependant, ce contre-exemple
fut obtenu de mani\`ere purement calculatoire et ne donne pas
beaucoup d'indications sur le probl\`eme. La question de la
dimension 5 reste toutefois ouverte (le cas des vari\'et\'es
toriques affines simpliciales de dimension 5 a \'et\'e r\'esolu
dans \cite{Si}). On montrera dans ce texte que la r\'eponse est
encore affirmative
pour les vari\'et\'es toriques de dimension 4.\\
\\
\indent On supposera fix\'e, pour toute la suite, $N \simeq
\zzz^n$ un r\'eseau de dimension $n$ (dans la derni\`ere section,
ce r\'eseau
sera de dimension 4) et $\Delta$ un \'eventail de $N \otimes \qqq$.\\
Introduisons maintenant les notations qui seront utilis\'ees ici.
On notera $\XC$ l'ensemble des points complexes ferm\'es de la
vari\'et\'e torique $X_{\Delta}$ associ\'ee \`a l'\'eventail
$\Delta$ et $\XR$ sa
partie r\'eelle.\\
On d\'esignera par $V = N \otimes \frac{\zzz}{2\zzz} =
\frac{N}{2N}$ la r\'eduction modulo 2 du r\'eseau $N$ et pour $p
\in \nnn$, $\Delta(p)$ d\'esignera le $p-$squelette de $\Delta$.
Explicitement:

$$\Delta(p) = \bigcup_{\sigma \in \Delta, dim(\sigma)=n-p} \ \sigma$$

\'Etant donn\'e un c\^one $\sigma$ de $\Delta$, on notera
$[\sigma]_2 \subset V$ la r\'eduction modulo 2 du sous r\'eseau
engendr\'e par $\sigma$ dans $N$. C'est \`a dire:

$$[\sigma]_2 = (Vect_{\qqq}(\sigma) \cap N) \otimes \zz$$

Si $\tau$ et $\sigma$ sont 2 c\^ones de $\Delta$, on notera
$p_{\tau,\sigma}:\Vtau \to \Vsig$ la projection naturelle et
$\pi_{\tau,\sigma}:H_0(\Vtau) \to H_0(\Vsig)$ l'application
induite par $p_{\tau,\sigma}$.\\

\section{Les suites spectrales r\'eelles et complexes}

L'homologie de Borel Moore \`a coefficients dans $\zz$ de la
vari\'et\'e torique r\'eelle $\XR$ est donn\'ee par le complexe
$A_*(\Delta)$ (cf \cite{Mc} pour les d\'etails) o\`u:

\begin{itemize}
    \item Pour $\sigma \in \Delta$, on pose $A_{\sigma}=
    H_0(\Vsig)$.

    \item Pour $0 \leq p \leq n$, on pose $A_p(\Delta)=
        \bigoplus_{\sigma \in \Delta(n-p)} A_{\sigma}$.
\end{itemize}

\noindent La diff\'erentielle $\partial_p:A_p(\Delta) \to
A_{p-1}(\Delta)$ est donn\'ee par la somme directe sur $\tau \in
\Delta(n-p-1)$ des applications:

$$ \sum_{\sigma > \tau,dim(\sigma)=n-p} \pi_{\tau,\sigma}:
H_0(\Vtau) \to \bigoplus_{\sigma > \tau,dim(\sigma)=n-p}
H_0(\Vsig)$$

Ce complexe correspond au complexe cellulaire obtenu par la
d\'ecomposition en orbite de $\XR$. Chaque orbite de dimension
$p$, $0 \leq p \leq n$, \'etant isomorphe \`a $(\rrr^*)^p$, est
une r\'eunion disjointe de $2^p$ quadrants de la forme $\rrr_{\leq
0}^p$ et ces derniers, tous hom\'eomorphes \`a $\rrr^p$,
constituent les
$p-$cellules de la d\'ecomposition cellulaire.\\
\\
\indent En ce qui concerne la vari\'et\'e torique complexe $\XC$,
la filtration donn\'ee par

$$\emptyset \subset X_0 \subset X_1 \subset \ldots \subset X_n =
\XC$$

\noindent o\`u pour $0 \leq p \leq n, X_p=X_\Delta(n-p)$ est la
r\'eunion des orbites de $\XC$ dont les dimensions sont
inf\'erieures o\`u \'egales \`a $p$ donne naissance \`a une suite
spectrale:
$$E^l_{pq}(\XC)\Rightarrow H_{p+q}(\XC)$$
convergeant vers l'homologie de Borel Moore de $\XC$.\\

Dans \cite{Jor}, A.Jordan montre que le terme $E^1_{pq}(\XC)$ est
donn\'e par:
 $$E^1_{pq}(\XC) \simeq \displaystyle \bigoplus_{\sigma \in
\Delta(n-p)} \Bq \Vsig$$ et que via ces identifications, la
diff\'erentielle d'ordre 1 de $E^1_{pq}(\XC)$, \\
$d^1_{pq}:E^1_{pq}(\XC) \to E^1_{p-1q}(\XC)$ est donn\'ee par la
somme directe sur $\Delta(n-p-1)$ des applications:

$$ \sum_{\sigma > \tau,dim(\sigma)=n-p} \Bq p_{\tau,\sigma}:
\Bq \Vtau \to \bigoplus_{\sigma > \tau,dim(\sigma)=n-p} \Bq
\Vsig)$$



Il est possible de d\'efinir une filtration $F_*$ (cf \cite{Mc}
pour les d\'etails) sur le complexe $A_*(\XR)$ de tel sorte que le
terme $\widetilde{G}^0_{pq}(\XR)$ de la suite spectrale
$\widetilde{G}^l_{pq}(\XR)$ associ\'ee \`a cette filtration soit
isomorphe au terme $E^1_{p+q,-p}(\XC)$. De plus, si l'on d\'esigne
par $\widetilde{\partial}^0_{pq}: \widetilde{G}^0_{pq}(\XR) \to
\widetilde{G}^0_{pq-1}(\XR)$ les diff\'erentielles du terme
$\widetilde{G}^0_{pq}(\XR)$, alors le diagramme suivant commute:

$$\xymatrix{
    \ar[d]^{d^1_{p+q,-p}} E^1_{p+q,-p}(X(\ccc)) \ar[r]^{\ \ \ \sim} &  \widetilde{G}^0_{pq}(X) \ar[d]_{\widetilde{\partial}^0_{pq}}\\
    E^1_{p+q-1,-p}(X(\ccc))\ar[r]^{\ \ \ \sim} & \widetilde{G}^0_{p,q-1}(X)
}$$

Ainsi, pour $p,q \in \zzz$, on obtient un isomorphisme:
\[ \widetilde{G}^1_{pq}(\XR) \simeq E^2_{p+q,-p}(\XC) \]
On poss\`ede ainsi un moyen de comparer $\sum_{p \in \nnn}
b_i(\XR)$ et $\sum_{p \in \nnn} b_i(\XC)$.

\begin{proposition}\label{smith}
Si la suite spectrale $\widetilde{G}^l_{pq}(\XR)$ d\'eg\'en\`ere
\`a l'ordre 1 (ie:
$\widetilde{G}^{\infty}_{pq}(\XR)=\widetilde{G}^1_{pq}(\XR)$),
Alors la suite spectrale $E^l_{pq}(\XC)$ d\'eg\'en\`ere \`a
l'ordre $2$ et on a en particulier:
$$\sum_{p \in \nnn} b_i(\XR) = \sum_{p \in \nnn} b_i(\XC)$$
\end{proposition}

\begin{preuve}
D'apr\`es l'in\'egalit\'e de Thom-Smith, on a:
$$ \sum_{p \in \nnn} b_i(\XR) \leq \sum_{p \in \nnn} b_i(\XC)$$
Si $\widetilde{G}^l_{pq}(\XR)$ d\'eg\'en\`ere \`a l'ordre 1, on a:
$$ \sum_{p \in \nnn} b_i(\XR) = \sum_{p,q} dim(\widetilde{G}^1_{pq}(\XR))$$
Comme $\sum_{p \in \nnn} b_i(\XC) = \sum_{p,q} dim(
E^{\infty}_{pq}(\XC))$, on obtient donc:
$$\sum_{p,q} dim(E^2_{pq}(\XC)) \leq \sum_{p,q}
dim(E^{\infty}_{pq}(\XC)) \ \ \Box$$
\end{preuve}

Il est montr\'e dans \cite{Mc} que la suite spectrale
$\widetilde{G}^l_{pq}(\XR)$ d\'eg\'en\`ere \`a l'ordre 1 si $\XR$
est de dimension inf\'erieure o\`u \'egale \`a 3 et dans
\cite{Ho3}, que c'est encore le cas
si $\XR$ est compacte \`a singularit\'es isol\'ees (en adaptant les r\'esultats de \cite{Br}).\\
On montrera dans la deuxi\`eme partie que c'est encore le cas
lorsque $\XR$ est de dimension 4.\\
\\
Comme les notations qui seront utilis\'ees par la suite
diff\`erent de celles choisies dans \cite{Mc}, rappelons
bri\`evement la
construction de la filtration $F_*$ sur $A_*(\XR)$.\\
\\
Pour construire la filtration $F_*$, il suffira de construire une
filtration sur chacun des espaces $H_0(\Vsig)$ pour $\sigma \in
\Delta$.\\

Soit $\sigma \in \Delta(n-p)$. Il est possible de munir
$H_0(\Vsig)$ d'une structure d'alg\`ebre par:

$$\begin{array}{crll}
  \times \colon & H_0(\Vsig)  \times H_0(\Vsig)  & \to & H_0(\Vsig)  \\
                & ([v_1],[v_2]) & \mapsto & [v_1 + v_2]
\end{array}$$

\begin{definition}
Si $H$ est un sous espace vectoriel de dimension $r$ de $\Vsig$,
on d\'efinit la classe de $H$ not\'ee $[H]$ dans $H_0(\Vsig)$ par:
$$[H]=\sum_{v \in H} \ [v]$$
Si $(e_i)_{i=1}^r$ est une base de $H$, on a \'egalement:
$$[H]=\prod_{i=1}^r ([0] + [e_i])$$
\end{definition}

Pour $0 \leq q \leq p$, on pose:
$$I^q=Vect([H], dim(H)=q)$$

\begin{remarque}
$I$ est le noyau du morphisme:

$$\begin{array}{crcl}
  \ \mu \colon & H_0(\Vsig) & \to & \displaystyle \zz \\
             & \sum_{i \in I} a_i [v_i]& \mapsto & \sum_{i \in I}
             a_i
\end{array}$$
En particulier $I$ est un id\'eal de $H_0(\Vsig)$ et $I^q$
correspond \`a la puissance $q^{\text{i\`eme}}$ de cet id\'eal.\\
\end{remarque}

Si $H$ est un sous espace vectoriel de dimension $q$ de $\Vsig$,
$q \in \nnn$, on notera $\overline{[H]}$ la classe de la cha\^ine
$[H]$ dans $\Iq H_0(\Vsig)$.\\
\\
La loi produit $\times$ sur $H_0(\Vsig)$ induit une structure
d'alg\`ebre gradu\'ee sur $Gr^*_I(H_0(\Vsig))$ d\'efinie sur les
g\'en\'erateurs par:

$$\begin{array}{crll}
  \times \colon & Gr^p_I(H_0(\Vsig)) \times Gr^q_I(H_0(\Vsig)) & \to & Gr^{p+q}_I(H_0(\Vsig))  \\
                & (\overline{[H]},\overline{[H']}) & \mapsto & \overline{[H + H']}
\end{array}$$

\noindent o\`u $H$ (resp $H'$) est un sous espace vectoriel de
$\Vsig$ de
dimension $p$ (resp $q$).\\
\\
Ensuite, il existe un isomorphisme d'alg\`ebre gradu\'ee:

$$ \psi_*: \BB \Vsig \to Gr^*_I(H_0(\Vsig)$$

\noindent On peut expliciter cet isomorphisme. Si on se donne $q
\in \nnn$, $1 \leq q \leq n$, alors $\psi_q$ est d\'efini sur les
g\'en\'erateurs par:

$$\begin{array}{crll}
  \psi_q  \colon & \Bq \Vsig  & \to &  Gr^q_I(H_0(\Vsig)) \\
                & v_1 \wedge\ldots \wedge v_q & \mapsto &
\overline{[H]}=\prod_{i=1}^q([0]+[v_i])
\end{array}$$

\noindent L'inverse de $\psi_q$ est alors donn\'e par:

$$\begin{array}{crll}
  \psi_q^{-1}  \colon & Gr^q_I(H_0(\Vsig)) & \to & \Bq \Vsig  \\
               & \overline{[H]} & \mapsto & v_1 \wedge\ldots \wedge v_q
\end{array}$$

\noindent o\`u $H$ est un sous espace vectoriel de $\Vsig$ de
dimension $q$ et $(v_1,\ldots,v_q)$ est une base de $H$.\\

Regardons par exemple ce qui se passe pour $q=1$.\\
Soit $v_1,v_2 \in V$. Alors:

   \[ ([0] + [v_1]) + ([0] + [v_2]) = ([0] + [v_1]) \times ([0] + [v_2])
   + ([0] + [v_1 + v_2]) \]

Ainsi si $\overline{([0] + [v_1])}$ d\'esigne la classe de $([0] +
[v_1])$ dans $\I H_0(\Vsig)$, $\overline{([0] + [v_2])}$ celle de
$([0] + [v_2])$ et $\overline{([0] + [v_1 + v_2])}$ celle de $([0]
+ [v_1 + v_2])$, alors:

\[ \overline{([0] + [v_1])} + \overline{([0] + [v_2])} = \overline{([0] + [v_1+v_2])} \]

et ainsi:

\[ \psi_1(v_1) + \psi_1(v_2) = \psi_1(v_1+v_2) \]

Finalement, la filtration $F_*$ sur le complexe $A_*(\XR)$ est
d\'efinie pour $q \leq 0$ par:

$$F_q(A_p(\XR))=\bigoplus_{\sigma \in \Delta(n-p)} I^{-q}
H_0(\Vsig)$$

Dans toute la suite, on consid\'erera le r\'earrangement
$G^l_{pq}(\XR)$ de la suite spectrale $\widetilde{G}^l_{pq}(\XR)$.
Ceci revient \`a poser:
\[ G^l_{pq}(\XR) = \widetilde{G}^l_{-q,p+q}(\XR)\]

\noindent Si l'on utilise cette notation, les diff\'erentielles
$\widetilde{\partial}^l_{pq}$ se transforment en:
$$\partial^l_{pq}:G^l_{pq}(\XR) \to G^l_{p-1,q+l}(\XR)$$

\section{\'Etude des diff\'erentielles du complexe $A_*(\XR)$}

On montrera \`a la fin de cette partie que la suite spectrale
$G^l_{pq}(\XR)$ d\'eg\'en\`ere \`a l'ordre 1 lorsque $\XR$ est une
vari\'et\'e torique r\'eelle de dimension $4$.

\subsection{M\'ethode}

En g\'en\'eral, pour montrer que la suite spectrale
$G^*_{pq}(\XR)$ d\'eg\'en\`ere \`a l'ordre 1, on utilisera le
proc\'ed\'e suivant:

\begin{definition}
On dira que le complexe $A_*(\XR)$ satisfait la propri\'et\'e
$s_{pq}$ pour $0 \leq p,q \leq n$ si la diff\'erentielle
$\partial_p$ v\'erifie:
$$ \forall c \in I^q(A_p(\XR)), \ \partial_p(c) \in I^{q+1}A_{p-1}(\XR), \ \exists d \in
I^{q+1}A_p(\XR), \ \partial_p(c) = \partial_p(d)$$
\end{definition}

Pour montrer que $A_*(\XR)$ satisfait la condition $s_{pq}$, on
utilisera la proposition suivante:

\begin{proposition}\label{prenul}
Pour $0 \leq p,q \leq n$, les 2 conditions suivantes sont
\'equivalentes:
\begin{itemize}
\item[(1)] Le complexe $A_*(\XR)$ satisfait la condition $s_{pq}$.

\item[(2)] Pour toute classe $\bar{c} \in G^0_{pq}(\XR) =
\displaystyle \Iq A_p(\XR)$ telle que
$\partial^0_{pq}(\overline{c})=0$, il est possible de trouver $c'
\in I^q A_p(\XR)$ v\'erifiant $\partial_p(c')=0$ et
$\overline{c}=\overline{c'} \in G^0_{pq}(\XR)$.
\end{itemize}
\end{proposition}

\begin{preuve}
Voir \cite{Si} ou \cite{Del2}.\\
\end{preuve}

\begin{remarque} \label{nultotal}
le complexe $A_*(\XR)$ satisfait la condition $s_{pq}$ si et
seulement si toutes les diff\'erentielles de la forme
$\partial_{pq}^r: G^r_{pq}(\XR) \to G^r_{p,q+r}(\XR)$ sont nulles
pour $r \geq 1$. Ceci revient \`a dire que l'on a $G^1_{pq}(\XR) =
G^{\infty}_{pq}(\XR)$.
\end{remarque}

\begin{corollaire}
Si pour tout $(p,q) \in \nnn^2$, le complexe $A_*(\XR)$ satisfait
la condition $s_{pq}$, alors la suite spectrale $G^*_{pq}(\XR)$
d\'eg\'en\`ere \`a l'ordre 1
($G^1_{**}(\XR)=G^{\infty}_{**}(\XR)$).
\end{corollaire}

Pour pouvoir appliquer la proposition \ref{prenul} et le
corollaire \ref{nultotal}, on aura besoin de bien comprendre la
relation entre diff\'erentielles de $G^0_{pq}(\XR)$ et celles du
complexe $A_*(\XR)$.

\subsection{Conditions $s_{p0}$, $p \in \nnn$, pour le complexe $A_*(\XR)$}

Soit $p \in \nnn$ fix\'e. On va montrer que le complexe $A_*(\XR)$
satisfait la condition $s_{p0}$. En fait ce r\'esultat est
d\'ej\`a d\'emontr\'e dans \cite{Mc} (remarque 7.4) avec une
d\'emonstration un peu diff\'erente mais la preuve de ce
r\'esultat que l'on va donner ici constitue un exemple simple
d'application de la m\'ethode qui sera utilis\'ee par la suite
pour montrer que $A_*(\XR)$ v\'erifie
d'autres conditions $s_{pq}$.\\

\begin{remarque}
Soit $c = (c_{\sigma})_{\sigma \in \Delta(n-p)} \in A_p(\XR)$.
Pour chaque $\sigma \in \Delta(n-p)$, il n'y a que 2
possibilit\'es pour la classe $\overline{c_{\sigma}}$ de
$c_{\sigma}$ dans $\displaystyle \frac{A_p(\XR)}{I}$:
\begin{itemize}

\item Soit $c_{\sigma} \in I$ et alors $\overline{c_{\sigma}}=0$.

\item Soit $c_{\sigma} \notin I$ et alors
$\overline{c_{\sigma}}=\overline{[0_{\sigma}]}$ o\`u $0_{\sigma}$
est l'\'el\'ement neutre de $\Vsig$ (dans ce cas, $[0_{\sigma}]$
est l'\'el\'ement unit\'e de $H_0(\Vsig)$).
\end{itemize}
\end{remarque}

\begin{lemme}\label{ligne1}
La condition $s_{p0}$ est v\'erifi\'ee par $A_*(\XR)$.
\end{lemme}

\begin{preuve}
Soit $c= (c_{\sigma})_{\sigma \in \Delta(n-p)} \in A_p(\XR)$, $c
\notin I(\A_p(\XR)$ tel que $\partial_p(c) \in I(A_{p-1}(\XR))$.\\
On note $\overline{c}= (\overline{c_{\sigma}})_{\sigma \in
\Delta(n-p)}$ la classe de $c$ dans $ \displaystyle
\frac{A_p(\XR)}{I}$. D'apr\`es la remarque pr\'ec\'edente, pour
tout $\sigma \in \Delta(n-p)$,
$\overline{c_{\sigma}}=\epsilon_{\sigma}
\overline{[0_{\sigma}]}$ avec $\epsilon_{\sigma} \in \zz$.\\
Soient $\gamma \in \Delta(n-p+1)$, $\sigma \in \Delta(n-p)$ tels
que $\gamma > \sigma$. Alors
$\pi_{\sigma,\gamma}([0_{\sigma}])=[0_{\gamma}]$ vu que
$\pi_{\sigma,\gamma}$ est un morphisme d'anneaux.\\
\\
Soit $\gamma \in \Delta(n-p+1)$. On consid\`ere les applications:
\begin{itemize}
    \item $p_{\gamma}: \displaystyle \bigoplus_{\sigma \in \Delta(n-p)}
    \frac{A_{\sigma}(\XR)}{I} \to \displaystyle \frac{A_{\gamma}(\XR)}{I}$
    d\'efinie par la somme directe sur\\ $\{\sigma \in \Delta(n-p), \gamma>\sigma
    \}$ des applications $p_{\sigma,\gamma}$.\\

    \item $\pi_{\gamma}:  \displaystyle \bigoplus_{\sigma \in \Delta(n-p)}
    A_{\sigma}(\XR) \to A_{\gamma}(\XR)$ d\'efinie par la somme directe sur\\
    $\{\sigma \in \Delta(n-p), \gamma>\sigma \}$ des applications $\pi_{\sigma,\gamma}$.\\
\end{itemize}

\noindent Avec ces notations, on a:

\begin{itemize}
    \item $p_{\gamma}=(\partial^1_{p0})_{|_{\frac{A_\sigma}{I}}}$

    \item $\pi_{\gamma}=(\partial_p)_{|_{A_{\sigma}}}$\\
\end{itemize}

Soit $c'=(c'_{\sigma})_{\sigma \in \Delta(n-p)} \in A_p(\XR)$
d\'efini par $c'_{\sigma}=\epsilon_{\sigma} [0_{\sigma}]$. Par
construction, on a donc $\overline{c} = \overline{c'} \in
\displaystyle \frac{A_p(\XR)}{I}$.\\
Comme $\partial_p(c) \in I(A_p(\XR))$,
$p_{\gamma}(\overline{c})=0$. En particulier, on doit avoir:

$$\sum_{\sigma \in \Delta(n-p), \sigma > \gamma}
\epsilon_{\sigma} = 0$$

\noindent Mais alors on doit aussi avoir:

$$\pi_{\gamma}(c')=0$$

\noindent En r\'esum\'e, $c'$ est un cycle et $\overline{c} =
\overline{c'}$ donc d'apr\`es \ref{prenul}, $s_{p0}$ est
v\'erifi\'ee.$\Box$
\end{preuve}

\subsection{Conditions $s_{nq}$,~ $q \in \nnn$, pour le complexe $A_*(\XR)$}

D'apr\`es la proposition \ref{prenul}, il suffira \`a nouveau,
pour $1 \leq q \leq n$ (la condition $s_{n0}$ est v\'erifi\'ee par
\ref{ligne1}), de montrer la proposition suivante: :

\begin{proposition} \label{precolonne1}
Soit $q \in \nnn^*$ fix\'e et $c \in I^q A_n(\XR)$ telle que
$\partial^0_{nq}(\overline{c})=0$ o\`u $\overline{c}$ d\'esigne la
classe de $c$ modulo $I^{q+1}$. Alors il existe $c'\in I^q
A_n(\XR)$ telle que $\partial_p(c')=0$ et
$\overline{c}=\overline{c'} \in G^0_{nq}(\XR) = \Iq A_n(\XR)$.
\end{proposition}

Il faut d'abord trouver un moyen naturel de relever un \'el\'ement
de $\displaystyle \frac{I^q}{I^{q+1}} A_n(\XR) \simeq \Bq V$ en un
\'el\'ement de $I^q A_n(\XR)$. Soit $(p,q) \in \nnn^2$ et $\sigma \in \Delta(p)$.\\

En g\'en\'eral, lorsque l'on fixe une base $e=(e_i)_{i=1 \ldots
n-p}$ de $\Vsig$, alors cela d\'efinit des applications
$s^q_{\sigma}:\Iq A_{\sigma} \to I^q(A_{\sigma})$, $0 \leq q \leq
n$, de la mani\`ere suivante:\\

Si $w \in \Bq \Vsig \simeq \Iq A_{\sigma}$, on peut d\'ecomposer
$w$ sous la forme:
$$w= \displaystyle \sum_{I \subset \{1 \ldots n-p\},
card(I)=q} \epsilon_I e_I$$

\noindent o\`u $I \subset \{1 \ldots n-p\}, \card(I)=q$ et $e_I=
\large \bigwedge_{i \in I} e_i$.\\
\\
Si $[e_I]$ d\'esigne la classe fondamentale du sous-espace de
dimension $q$ de $V$ engendr\'e par les vecteurs $e_i$, $i \in I$
dans $V$ (ie: $[e_I]=\displaystyle \sum_{v \in Vect(e_i,i \in I)}
[v])$, alors on pose:
$$s^q_{\sigma}(w)=\sum_{I \subset \{1 \ldots n-k\},card(I)=q}
\epsilon_I [e_I]$$

\begin{remarque}
L'application $s^q_{\sigma}$ est une section de la projection
$f^q_{\sigma} :I^q(A_{\sigma})\to \Iq A_{\sigma}$ ($f^q_{\sigma}
\circ s^q_{\sigma} = id$). Par contre, $s^q_{\sigma} \circ
f^q_{\sigma} \neq id$. Pour le voir, on peut prendre par exemple
$\sigma=\{0\}$, $q=2$ et $n=3$. Soit $(e_i)_{i=1\ldots 3}$ une
base de $V$, $s^2_{\{0\}}:\displaystyle \frac{I^2}{I^3} H_0(V) \to
I^2 H_0(V)$ l'application d\'efinie par cette base et
$f^2_{\{0\}}:I^2 H_0(V) \to \displaystyle \frac{I^2}{I^3} H_0(V)$
la projection. Si $c= ([0] + [e_1+e_2])([0]+[e_3])$ alors:

\begin{eqnarray*}
s^2_{\{0\}} \circ f^2_{\{0\}}(c) & = & s^2_{\{0\}}( (e_1+e_2) \wedge e_3) \\
              & = & s^2_{\{0\}}(e_1 \wedge e_3 + e_2 \wedge e_3)\\
              & = & ([0]+[e_1])([0]+[e_3]) + ([0]+[e_2])([0]+[e_3])\\
              & \neq & c
\end{eqnarray*}
\end{remarque}

\begin{remarque}
Lorsque $q=1$, il est possible de d\'efinir une telle application
sans fixer au pr\'ealable une base de $V$.\\
En effet, si $v \in \displaystyle \frac{V}{[\sigma]_2} \simeq
\frac{I}{I^2}~ A_{\sigma}$, il suffit de poser
$s^1_{\sigma}(v)=[0]+[v]$.\\
\end{remarque}

Il faut ensuite s'assurer que la base $e$ que l'on choisi soit
bonne. Plus pr\'ecisemment, on veut que l'application $s_e$
v\'erifie la propri\'et\'e suivante:

$$ \partial^0_{nq}(\overline{c})=0 \in \Iq A_n(\XR) \Rightarrow
\partial_n(s_e(\overline{c}))=0 \in I^q A_n(\XR)$$

On aura besoin pour cela du lemme suivant:

\begin{lemme} \label{annulation}
Soit $[H]$ la classe fondamentale d'un hyperplan de $V$ dans
$H_0(V)$ et soit $V'$ un sous-espace vectoriel de $V$.\\
On consid\`ere l'application $\pi: H_0(V) \to H_0(\frac{V}{V'})$
induite par la projection naturelle $p: V \to \frac{V}{V'}$ de $V$
sur $V'$. Alors:
$$H \cap V' \neq \{0 \} \Rightarrow \pi([H])=0$$
\end{lemme}

\begin{preuve}
Supposons que $H \cap V'$ contienne $v_1 \neq 0$. On peut
compl\'eter $v_1$ en une base $(v_1, \ldots, v_m)$ de $H$ o\`u
$m=\dimm(H)$. Alors on a dans $H_0(V)$:
$$[H] = \prod_{i=1}^m ([0] + [v_i])$$
Ainsi:
\begin{eqnarray*}
    \pi([H]) & = & \prod_{i=1}^m (\pi([0]) + \pi([v_i]))\\
             & = & \prod_{i=1}^m ([0] + [p(v_i)])\\
             & = & ([0] + [0]) \times \prod_{i=2}^m ([0] +
             [p(v_i)]) \\
             & = &  0 ~\Box
    \end{eqnarray*}
\end{preuve}

\noindent Toute la difficult\'e sera donc d'exhiber une telle base.\\
\\
\indent On pose $\Delta(1)=\bigcup_{i=1}^m \tau_i$ et pour chaque
c\^one $\tau_i$, on notera $v_i$ le g\'en\'erateur de
$[\tau_i]_2$.
D\'emontrons \`a pr\'esent la proposition \ref{precolonne1}.\\

\begin{preuve}
\indent Comme $\partial_n(c) \in I^{q+1}(A_{n-1}(\Delta))$, on
doit avoir:

$$\begin{array}{crcl}
   \partial^0_{nq} \colon & G^0_{nq}(\XR) \simeq \bigwedge^q V & \to & G^0_{n-1,q}(\XR) \simeq \bigoplus_{i=1}^r \bigwedge^q \frac{V}{[\tau_i]_2} \\
       &  \overline{c} & \mapsto & 0
\end{array}$$

\noindent On peut supposer que $v_1, \ldots, v_r$ est une base de
$\Vect(v_i,
\ i=1, \ldots, m)$.\\
Il y a deux cas:
\begin{itemize}

    \item[*] Si $r>q$, alors $\partial^0_{nq}$ est injective.
    Ainsi $\overline{c}=0$ et $c \in I^{q+1}$. On peut donc
    choisir $c'=0$.

    \item[*] Si $r \leq q$, on choisit une base $e=(e_1, \ldots,
    e_n)$ de $V$ telle que $e_i=v_i$ pour $i=1, \ldots,r$.\\

    On peut alors \'ecrire $\overline{c}$ dans cette base sous la forme:

    $$\overline{c}=\sum_I \epsilon_I e_I$$

    Comme $\partial^0_{nq}(\overline{c})=0$, chacun des termes $e_I$ doit v\'erifier:

    $$\epsilon_I \neq 0 \Rightarrow \{1,\ldots, r\} \subset I$$

    Maintenant qu'une base de $V$ est fix\'ee, on dispose de
    l'application $s_e$ qui lui est associ\'ee. On pose
    $c'=s_e(\overline{c}) \in I^q(A_n(\Delta))$.\\

    Si $\overline{c'} \in G^0_{nq}$ d\'esigne la classe de $c'$ modulo $I^{q+1}$, on a ainsi:
    $$\overline{c} = \overline{c'}$$

    Il reste \`a montrer que $\partial_n(c')=0$. Or pour
    $i=1,\ldots,m$ et $I \subset \{1, \ldots, m\}$, $\card(I)=q$, v\'erifiant $\epsilon_I \neq 0$,
    le sous-espace $\Vect(e_j, \ j \in I)$ contient $v_i$.\\
    Ainsi d'apr\`es le lemme \ref{annulation}, on a pour $i=1,\ldots,m$:

    $$\begin{array}{crcl}
    \pi_i \colon & H_0(V) & \to & H_0(\frac{V}{[\tau_i]_2}) \\
                    & [e_I] & \mapsto & 0
    \end{array}$$

    Alors on a:
    \begin{eqnarray*}
    \pi_i(c') & = & \pi( \sum_I \epsilon_I [e_I]) \\
              & = & \sum_I \epsilon_I \pi_i([e_I])\\
              & = & 0
    \end{eqnarray*}

    Ainsi $\partial_n(c')=0$ et $c'$ convient. $\Box$
\end{itemize}
\end{preuve}

\subsection{Condition $s_{n-1,1}$ pour le complexe $A_*(\XR)$}

Toujours d'apr\`es \ref{prenul}, il faudra montrer la proposition
suivante:

\begin{proposition}\label{fleche}
Soit $c \in I A_{n-1}(\XR)$ telle que
$\partial^0_{n-1,1}(\overline{c})=0$ o\`u $\overline{c}$ d\'esigne
la classe de $c$ modulo $I^2$. Alors il existe $c'\in I^2
A_{n-1}(\XR)$ telle que $\partial_p(c')=0$ et
$\overline{c}=\overline{c'} \in G^0_{n-1,1}(\XR) = \I A_n(\XR)$
\end{proposition}

Rappelons que la diff\'erentielle \`a laquelle on s'int\'eresse
est:

$$\partial^0_{n-1,1}:G^0_{n-1,1}(\XR) \simeq \bigoplus_{\tau \in \Delta(1)} \Vtau  \to G^0_{n-2,1}(\XR)
\simeq \bigoplus_{\sigma \in \Delta(2)} \Vsig $$

\begin{preuve}
Soit $\Delta(1)=\bigcup_{i=1}^m \tau_i$. Pour $1 \leq i \leq m$,
on notera $\rho_i \in N$ le vecteur primitif engendrant le c\^one
$\tau_i$ et $v_i \in V= \displaystyle \frac{N}{2N}$ d\'esignera sa classe modulo 2.\\
\\
\'Etant donn\'e $\sigma \in \Delta(2)$, il existe $1 \leq i<j \leq
m$ tels que $\sigma$ soit engendr\'e par $\rho_i$ et $\rho_j$. On
pose
alors $\sigma=\sigma_{ij}$.\\
\\
Soit $c \in I(A_{n-1}(\XR))$ telle que $\partial_{n-1}(c) \in
\displaystyle \frac{I}{I^2} A_{n-2}(\XR)$. On pose $c=(c_i)_{i=1
\ldots m}$ et on d\'esigne par $\overline{c_i} \in \displaystyle
\frac{V}{[\tau_i]_2}$
la classe de $c_i$ modulo $I^2$.\\
\\
Pour chaque $i=1 \ldots m$, on d\'esigne par $s_i$ l'application
d\'efinie par:

$$\begin{array}{crcl}
  s_i \colon & \displaystyle \frac{V}{[\tau_i]_2} & \to & \displaystyle I(H_0(\frac{V}{[\tau_i]_2})) \\
       & v & \mapsto & [0] + [v]
\end{array}$$

\noindent et l'on note $c_i'=s_i(c_i)$.\\
On obtient ainsi une nouvelle cha\^ine $c' \in A_{n-1}(\XR) $ qui
v\'erifie:

\begin{itemize}
    \item $c' \in I(A_{n-1}(\XR))$.
    \item $\overline{c}=\overline{c'}$ modulo $I^2$ si
    $\overline(c')$ d\'esigne la classe de $c'$ modulo $I^2$.\\
\end{itemize}

Il reste \`a v\'erifier que $c'$ v\'erifie $\partial_{n-1}(c')=0$.
Soit $b=\partial_{n-1}(c')$. On pose $b=(b_{\sigma})_{\sigma \in
\Delta(2)}$.\\
\'Etant donn\'e $\sigma_{ij} \in \Delta(2)$, on notera:

\begin{itemize}
    \item $p_{ij}: \displaystyle \frac{V}{[\tau_i]_2} \to \displaystyle \frac{V}{[\sigma_{ij}]_2}$
    la projection naturelle.

    \item $\pi_{ij} :\displaystyle H_0(\frac{V}{[\tau_i]_2}) \to
    \displaystyle H_0(\frac{V}{[\sigma_{ij}]_2})$ l'application induite par
    $p_{ij}$.
\end{itemize}

Soit $\sigma = \sigma_{ij} \in \Delta(2)$. On a:
\begin{eqnarray*}
b_{\sigma_{ij}} & = & \pi_{ij}(c_i') + \pi_{ji}(c_j') \\
                & = & \pi_{ij}([0] + [\overline{c_i}]) + \pi_{ji}([0] + [\overline{c_j}])\\
                & = & [0] + [p_{ij}(\overline{c_i})] + [0] + [p_{ji}(\overline{c_j})]\\
\end{eqnarray*}

Comme $\partial_{n-1}(c) \in I^2$, on doit avoir
$p_{ij}(\overline{c_i}) + p_{ji}(\overline{c_j})=0$ donc
$b_{\sigma_{ij}}=0$ et finalement, on obtient bien $b=0$.$\Box$
\end{preuve}

\subsection{Maximalit\'e des vari\'et\'es toriques de dimension 4}

\begin{theoreme}
Soit $X$ une vari\'et\'e torique de dimension $4$. Alors la suite
spectrale $G^l_{pq}(\XR)$ associ\'ee \`a $X$ d\'eg\'en\`ere \`a
l'ordre 1.
\end{theoreme}

\begin{preuve}

\noindent D'apr\`es \ref{ligne1}, les conditions $s_{00}, s_{10},
s_{20}, s_{30}$ et $s_{40}$ sont v\'erifi\'ees par $A_*(\XR)$. la
proposition \ref{precolonne1} permet de voir que les conditions
$s_{41}, s_{42}, s_{43}, s_{44}$ le sont \'egalement. Ensuite, la
proposition \ref{fleche} donne la condition $s_{31}$. Finalement,
les autres conditions $s_{pq}$ \'etant trivialement v\'erifi\'ees
lorsque $n=4$, toutes les diff\'erentielles $\partial^*_{pq}$ de
la suite spectrale $G^*_{pq}(\XR)$ de degr\'e sup\'erieur \`a 1
sont nulles. Ainsi, $G^*_{pq}(\XR)$ d\'eg\'en\`ere \`a l'ordre 1
et d'apr\`es \ref{smith}, $\XC$ est maximale.
\end{preuve}


\begin{thebibliography}{99}

\bibitem{Br} M.Brion, {\em The structure of the polytope
algebra}, T\^ohoku Math. J. 49 (1997), 1-32.

\bibitem{Mc} F. Bihan, M. Franz, C. Mc Crory, J. Van Hamel, {\em Is every toric variety
    an M-variety?}, Manuscripta Math 120 (2006), 217-232.

\bibitem{Del2} P. Deligne, {\em Th\'eorie de Hodge II}, Publ. Math.
IHES 40 pp 5-58 (1971).

\bibitem{Ful} W. Fulton, {\em Introduction to toric variety},
    Princeton University Press, Princeton, NJ, 1993.

\bibitem{Ho} V. Hower, {\em Hodge spaces of real toric varieties},
arXiv: alg-geom/0705.0516, 2007.

\bibitem{Ho2} V. Hower, {\em A counterexample to the
maximality of toric varieties}, arXiv:alg-geom/0611925, 2006.

\bibitem{Ho3} V.Hower, {\em Hodge spaces of real toric varieties},
Ph .D thesis, University of Georgia, 2007.

\bibitem{Jor} A. Jordan, {\em Homology and Cohomology of toric
varieties}, Ph. D. thesis, Universit\"at Konstanz, 1997.

\bibitem{Si} A. Sine, {\em Probl\`eme de maximalit\'e pour les vari\'et\'es toriques},
Ph. D. thesis, Universit\'e d'Angers, 2007.

\end{thebibliography}
\end{document}